\documentclass[12pt]{article}

\usepackage{pst-poly}     
\usepackage{pst-all}      
\usepackage{amsfonts,amssymb,amsmath,latexsym,bm}

\newtheorem{thm}{Theorem}[section]

\newtheorem{pps}[thm]{Proposition}

\newtheorem{lem}[thm]{Lemma}
\newtheorem{cjt}[thm]{Conjecture}

\newenvironment{pf}[1][Proof]{\noindent\textbf{#1.} }{\hfill\rule{1mm}{2mm}}

\makeatletter \@addtoreset{equation}{section} \makeatother

\begin{document}

\title{Pebbling on $C_{4k+3}\times G$ and $M(C_{2n})\times G$
\thanks{Supported by ``the Fundamental Research Funds for the Central Universities" and the NSF of the People's Republic of China(Grant
 No. 61272008, No. 11271348 and
 No. 10871189).}
  }
\author
{Zheng-Jiang Xia,\quad Yong-Liang Pan\footnote{Corresponding
author:
ylpan@ustc.edu.cn},\quad Jun-Ming Xu\ \\ \\
{\small     School of Mathematical Sciences,} \\
{\small             University of Science and Technology of China,}   \\
{\small             Hefei, Anhui, 230026, P. R. China}\\
{\small             Email: ylpan@ustc.edu.cn}  \\
}

\date{}
\maketitle

\begin{quotation}
\noindent\textbf{Abstract}: The  pebbling  number  of  a  graph  $G$, $f(G)$,  is  the  least  $p$  such  that,  however  $p$  pebbles  are  placed  on  the  vertices  of  $G$,  we  can  move  a  pebble  to  any  vertex by  a  sequence  of  moves, each  move  taking  two  pebbles  off one  vertex  and  placing one  on  an  adjacent  vertex. It is conjectured that for all graphs $G$ and $H$, $f(G\times H)\leq f(G)f(H)$.
If the graph $G$ satisfies the odd two-pebbling property, we will prove that $f(C_{4k+3}\times G)\leq f(C_{4k+3})f(G)$ and $f(M(C_{2n})\times G)\leq f(M(C_{2n}))f(G)$, where $C_{4k+3}$ is the odd cycle of order $4k+3$  and $M(C_{2n})$
is the middle graph of the even cycle $C_{2n}$.\\


\end{quotation}

\section{Introduction}

Pebbling in graphs was first introduced by Chung~\cite{c89}. Consider a connected graph with a fixed number of pebbles
distributed on its vertices. A pebbling move consists of the removal of two pebbles from a vertex and the placement of
one pebble on an adjacent vertex. The pebbling number of a vertex $v$, the target vertex, in a graph $G$ is the smallest number $f(G,v)$
with the property that, from every placement of $f(G,v)$ pebbles on $G$, it is possible to move one pebble to $v$ by a sequence of
pebbling moves. The \emph{t-pebbling number}
of $v$ in $G$ is defined as the smallest number $f_t(G,v)$ such that from every placement of $f_t(G,v)$ pebbles, it is possible to move
$t$ pebbles to $v$. Then the \emph{pebbling number} and the \emph{t-pebbling number} of $G$ are the smallest numbers, $f(G)$ and
$f_t(G)$, such that from any placement of $f(G)$ pebbles or $f_t(G)$ pebbles, respectively, it is possible to move one or $t$
pebbles, respectively, to any specified target vertex by a sequence of pebbling moves. Thus, $f(G)$ and $f_t(G)$ are the
maximum values of $f(G,v)$ and $f_t(G,v)$ over all vertices $v$.\\


Chung~\cite{c89} defined the two-pebbling property of a graph, and Wang~\cite{w01} extended her definition to the odd two-pebbling property as follows.\\

\indent Suppose $p$ pebbles are located on $G$, let $k$ be the number of occupied vertices (vertices with at least one pebble),
and $r$ be the number of vertices with an odd number of pebbles. $G$ satisfies the \emph{two-pebbling property
} means two pebbles can be moved to any vertex of $G$ whenever $p>2f(G)-k$, and the \emph{odd two-pebbling property}
means two pebbles can be moved to any vertex of $G$ whenever $p>2f(G)-r$. It is clear that the graph which satisfies the
two-pebbling property also satisfies the odd two-pebbling property, and there exist graphs do not satisfy the two-pebbling property, which are called the Lemke graphs.

\indent The {\it middle graph} of a graph $G$, denoted by $M(G)$, is obtained from $G$ by inserting a new vertex into each edge of $G$,
and joining the new vertices by an edge if the two edges they inserted share the same vertex of $G$.
For any two graphs $G$ and $H$, we define the \emph{Cartesian product} $G\times H$ to be the graph with vertex
set $V(G\times H)$ and edge set the union of $\{((a, v),  (b, v))| (a, b)\in E(G),  v\in E(H)\}$  and $\{((u, x),  (u, y))|u\in V(G),  (x, y)\in E(H)\}$.

The following conjecture~\cite{c89}, by Ronald Graham, suggests a constraint on the pebbling number
of the product of two graphs.

\begin{cjt}[Graham]
The pebbling number of $G\times H$ satisfies $f(G\times H)\leq f(G)f(H)$.
\end{cjt}

There are a few results verify Graham's conjecture. It has been proved that Graham's conjecture holds for a tree by a graph with the odd-two-pebbling property, especially a tree by a tree~\cite{m92}, an even cycle by a graph with the odd two-pebbling property, a cycle by a cycle~\cite{h03}, a complete graph by a graph with the two-pebbling property~\cite{c89} and a complete bipartite graph by a graph with the two-pebbling property~\cite{fk01}, a fan graph by a fan graph and a wheel graph by a wheel graph~\cite{fk02}, a thorn graph of the complete graph with every $p_i>1$ by a graph with the two-pebbling property~\cite{wz09}, and the middle graph of an odd cycle by the middle graph of a cycle~\cite{yzz}.\\

In Section $2$, we show  that Graham's conjecture holds for the odd cycle $C_{4k+3}$ by a graph with the odd
two-pebbling property.\\

In Section $3$, we show that Graham's conjecture holds for the middle graph of an even cycle by a graph with the odd two-pebbling property, especially, the middle graph of an even cycle by the middle graph of an even cycle.\\

Given a distribution of pebbles on $G$, Let $p(K)$  be the number of pebbles on a subgraph $K$ of $G$, $p(v)$  be the number of pebbles on
vertex $v$ of $G$ and $k(K)$ ($r(K)$) to be the number of vertices of $K$ with at least one pebble (with an odd number of pebbles). Moreover, denote by $\tilde{p}(K)$ and $\tilde{p}(v)$ the number of pebbles on $K$ and $v$ after some sequence of pebbling moves, respectively.

Let $T$ be a tree with a specified vertex $v$. $T$ can be viewed as a directed tree denoted
by $\vec{T}_v$ with edges directed toward a specified vertex $v$, also called the root. A path-partition
is a set of non-overlapping directed paths the union of which is $\vec{T}_v$.
The  path-size  sequence  of  a  path-partition  ${P_1,\ldots, P_n,}$  is  an  $n-$tuple $(a_1,\ldots, a_n)$,  where  $a_j$  is
the  length  of  $P_i$ (i.e.,  the  number  of  edges  in  it), with $a_1\geq a_2\geq\ldots\geq a_n$.
A path-partition is said to majorize another if the nonincreasing sequence
of its path size majorizes that of the other. That is,
$(a_1,a_2,\ldots,a_r)>(b_1,b_2,\ldots,b_t)$ if and only if $a_i > b_i$ where
$i= \min\{j:a_j\neq b_j \}$.
A path-partition of a tree T is said to be maximum if it majorizes all other
path-partitions.\\

The following two Lemmas will be the key tools in the next sections.
\begin{lem}{\rm{(\cite{c89})}}\label{lem1}
The pebbling number $f_t(T,v)$ for a vertex $v$ in a tree $T$ is
$t2^{a_1}+2^{a_2}+\cdots+2^{a_r}-r+1$, where $a_1, a_2,\ldots,a_r$ is the sequence of the path sizes in a
maximum path-partition of $\vec{T}_v$.
\end{lem}

\begin{lem}{\rm{(\cite{m92})}}\label{lem2.0}
If $T$ is a tree, and $G$ satisfies the odd two-pebbling property, then $f(T\times G,(x,g))\leq f(T,x)f(G)$ for every vertex $g$ in $G$.
In particular, if $P_m=x_1x_2\ldots x_m$ be a path, then $f(P_m\times G,(x_i,g))\leq f(P_m,x_i)f(G)
=(2^{i-1}+2^{m-i}-1)f(G)\leq2^{m-1}f(G)$.
\end{lem}

\section{$C_{4k+3}\times G$}

In 2003,  D. S. Herscovici proved the following two theorems about cycles.
\begin{thm}\rm{(\cite{h03})}
If $G$ satisfies the odd two-pebbling property, then
$$f(C_{2n}\times G)\leq f(C_{2n})f(G)=2^nf(G).$$
\end{thm}

\begin{thm}\rm{(\cite{h03})}
Suppose $G$ is a graph with $m\geq 5$ vertices which satisfies the odd two-pebbling property and the following inequality
\begin{align}\label{eq1.0}
4f_4(G)<14f(G)-2(m-5),
\end{align}
then, $f(C_{2n+1}\times G)\leq f(C_{2n+1})f(G)$ for $n\geq3$.
\end{thm}

The inequality (\ref{eq1.0}) holds for all odd cycles, but does not hold even for paths or even cycles.
In this section, we show the following theorem.
\begin{thm}\label{thm1.1}
If $G$ satisfies the odd two-pebbling property, then
$$f(C_{4k+3}\times G)\leq f(C_{4k+3})f(G).$$
\end{thm}

Throughout this section, we use the following notation.

Let  the vertices of $C_{4k+3}$ be $\{v_0,v_1,\ldots,v_{4k+1},v_{4k+2}\}$ in order.
We define the vertex subsets $A$ and $B$ of $C_{4k+3}$ by
 $$A=\{v_1,v_2,\ldots,v_{2k}\},B=\{v_{2k+3},v_{2k+4},\ldots,v_{4k+2}\}$$

For simplicity, among $C_{4k+3}\times G$, let $p_i=p(v_i\times G)$, $r_i=r(v_i\times G)$, $p(A)=p(A\times G)$, $p(B)=p(B\times G)$. Thus, the number of pebbles in a
distribution on $C_{4k+3}\times G$ is given by $p_0+p(A)+p(B)+p_{2k+1}+p_{2k+2}$.

\begin{lem}{\rm{(\cite{psv95})}}
The pebbling numbers of the odd cycles $C_{4k+1}$ and $C_{4k+3}$ are
\begin{align*}
f(C_{4k+1})=&\frac{2^{2k+2}-1}{3}
=1+2^2+2^4+\cdots+2^{2k}.\\
f(C_{4k+3})=&\frac{2^{2k+3}+1}{3}
=1+2^1+2^3+\cdots+2^{2k+1}.
\end{align*}
\end{lem}

\begin{lem}\label{lem1.1}
Let $P_{2k}=x_1x_2\ldots x_{2k}$ be a path with length $2k-1$, and let $g$ be some vertex in a graph $G$ which satisfies the odd two-pebbling property. Then, from any arrangement of $(2^1+2^3+\cdots+2^{2k-1})f(G)$ pebbles on $P_{2k}\times G$, it is possible to put a pebble on every $(x_i,g)$ at once, where $i=1,3,\ldots,2k-1$.
\end{lem}
\begin{pf}
We use induction on $k$, where the case $k=1$ is trivial.

Suppose that there are $(2^1+2^3+\cdots+2^{2k-1})f(G)$ pebbles on $P_{2k}\times G$. Then there are at least $(2^1+2^3+\cdots+2^{2k-3})f(G)$ pebbles on $\{x_3,x_4,\ldots,x_{2k}\}\times G$
(or on $\{x_1,x_2,\ldots,x_{2k-2}\}\times G$). By induction, we can use these pebbles to put one pebble to each of these vertices $\{(x_3,g),(x_5,g),\ldots,(x_{2k-1},g)\}$(or $\{(x_1,g),(x_3,g),\ldots,(x_{2k-3},g)\}$). By Lemma~\ref{lem1}, $f(P_{2k},x_1)=2^{2k-1}$, $f(P_{2k},x_{2k-1})=2^{2k-2}+1\leq2^{2k-1}$.
Thus with the remaining $2^{2k-1}f(G)$ pebbles, one pebble can be moved to $(x_1,g)$(or $(x_{2k-1},g)$), and we are done.
\end{pf}\\

Similarly, we can get the following Lemma.
\begin{lem}\label{lem1.2}
Let $P_{2k+1}=x_1x_2\ldots x_{2k+1}$ be a path with length $2k$, and let $g$ be some vertex in a graph $G$ which satisfies the odd two-pebbling property. Then, from any arrangement of $(2^2+2^4+\cdots+2^{2k})f(G)$ pebbles on $P_{2k+1}\times G$, it is possible to put a pebble on every $(x_i,g)$ at once, where $i=1,3,\ldots,2k-1$.
\end{lem}

From the proof of Theorem 3.2 in [5], it follows that
\begin{lem}{\rm{(\cite{h03})}}\label{lem1.3}
If $p(A)\geq2^{2k-1}f(G)$, then with $f(C_{4k+3})f(G)$ pebbles on $C_{4k+3}\times G$, one pebble can be moved to $(v_0,g)$.
\end{lem}

\noindent{\bf\emph{Proof of Theorem 2.3}:}\\
Suppose that there are $f(C_{4k+3})f(G)$ pebbles located on $C_{4k+3}\times G$, then
\begin{align}\label{eq1.1}
p_0+p_{2k+1}+p_{2k+2}+p(A)+p(B)=(1+2^1+2^3+\cdots+2^{2k+1})f(G).
\end{align}

Without loss of generality, we may assume that $p(A)\geq p(B)$ and the target vertex is $(v_0,g)$. The case $k=0$ is trivial, so we  assume that $k\geq1$.

Note that the vertices of $B\cup\{v_{2k+2}\}\cup\{v_0\}$ form a path isomorphic to $P_{2k+2}$. It follows from Lemma 1.3 that if we move as many as possible pebbles from $v_{2k+1}\times G$ to $v_{2k+2}\times G$, then one pebble could be moved to $(v_0,g)$ unless
\begin{align}\label{eq1.2}
\frac{p_{2k+1}-r_{2k+1}}{2}+p_{2k+2}+p(B)+p_0< 2^{2k+1}f(G).
\end{align}

From Lemma~\ref{lem1.3}, we could move one pebble to $(v_0,g)$ unless
\begin{align}\label{eq1.3}
p(A)< 2^{2k-1}f(G).
\end{align}

If (\ref{eq1.2}) holds, then
\begin{align}\label{eq1.4}
\frac{p_{2k+1}+r_{2k+1}}{2}+p(A)>(1+2^1+2^3+\cdots+2^{2k-1})f(G).
\end{align}

From (\ref{eq1.3}) and (\ref{eq1.4}), we can get
$p_{2k+1}+r_{2k+1}>2f(G),$
and
\begin{align*}
\frac{p_{2k+1}-(2f(G)-r_{2k+1}+2)}{2}+p(A)\geq (2^1+2^3+\cdots+2^{2k-1})f(G).
\end{align*}

 This implies that we can move enough pebbles from $v_{2k+1}\times G$ to $A\times G$ so that the number of the pebbles on $A\times G$ will reach to $(2^1+2^3+\cdots+2^{2k-1})f(G)$, and at the same time there are $h_{2k+1}$  pebbles are kept on $v_{2k+1}\times G$, where
 $$h_{2k+1}=\begin{cases}2f(G)-r_{2k+1}+2, & \mbox{if}\;\;  r_{2k+1}\geq 2,\\
 2f(G),& \mbox{if}\;\;  r_{2k+1}\leq 1.\end{cases}$$
Assume that $2x$ pebbles are taken away from $v_{2k+1}\times G$ such that there are $x$ pebbles reaching $A\times G$, i.e.,
\begin{align}\label{eq1.5}
x+p(A)=(2^1+2^3+\cdots+2^{2k-1})f(G).
\end{align}

\indent   Step 1.\; With the $h_{2k+1}$ pebbles  on $v_{2k+1}\times G$, we can move one pebble to $(v_{2k},g)$. \\

Now  there are $p_{2k+1}-2x- h_{2k+1}$ pebbles on $v_{2k+1}\times G$, namely,
\begin{align*}
\tilde{p}_{2k+1}=&p_{2k+1}-2x-h_{2k+1}\\
=&p_{2k+1}+2p(A)-(2^2+\cdots+2^{2k})f(G)-h_{2k+1}.
\end{align*}

\noindent So the remaining pebbles on $\{v_0,v_{2k+1},v_{2k+2}\}\times G$ is
\begin{align*}\label{eq1.6}
&p_0+p_{2k+2}+\tilde{p}_{2k+1}\\
&=p_0+p_{2k+2}+p_{2k+1}+2p(A)-(2^2+\cdots+2^{2k})f(G)-h_{2k+1}\\
&\geq p_0+p_{2k+2}+p_{2k+1}+p(A)+p(B)-(2^2+\cdots+2^{2k})f(G)-h_{2k+1}\\
&\geq(1+2^2+2^4+\cdots+2^{2k})f(G).
\end{align*}


For $p_0<f(G)$ (otherwise one pebble can be moved to $(v_0,g)$, and we are done), so
\begin{align}
p_{2k+2}+\tilde{p}_{2k+1}\geq(2^2+2^4+\cdots+2^{2k})f(G).
\end{align}

\indent Step 2.\, It follows from (\ref{eq1.5}) and Lemma~\ref{lem1.1} that with $(2^1+2^3+\cdots+2^{2k-1})f(G)$ pebbles on $A\times G$, we can put one pebble to each vertex of $\{(v_1,g),(v_3,g),\ldots,(v_{2k-1},g)\}$.\\
\indent Step 3.\, From the inequality (\ref{eq1.6}) and Lemma~\ref{lem1.2}, it follows that, with $(2^2+2^4+\cdots+2^{2k})f(G)$ pebbles on $\{v_{2k+1},v_{2k+2}\}\times G$, we can put one pebble to each vertex of $\{(v_2,g),(v_4,g),\ldots,(v_{2k},g)\}$.\\
\indent  The above three steps implies that  at least one pebble can be moved to $(v_0,g)$. \hfill\rule{1mm}{2mm}

\section{$M(C_{2n})\times G$}

%

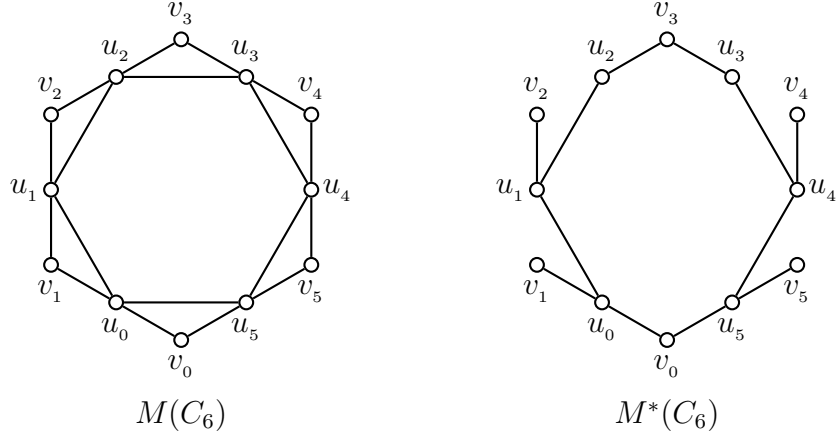
\begin{figure}[ht]
\begin{center}
\hspace*{30pt}
\begin{pspicture}(-5,-1)(5,5.5)

\cnode(-3.73,4.5){3pt}{x1}\cnode(-2.865,4){3pt}{x2}\cnode(-2,3.5){3pt}{x3}\cnode(-2,2.5){3pt}{x4}\cnode(-2,1.5){3pt}{x5}
\cnode(-2.865,1){3pt}{x6}\cnode(-3.73,0.5){3pt}{x7}\cnode(-4.595,1){3pt}{x8}\cnode(-5.46,1.5){3pt}{x9}\cnode(-5.46,2.5){3pt}{x10}
\cnode(-5.46,3.5){3pt}{x11}\cnode(-4.595,4){3pt}{x12}
\ncline{x1}{x2}\ncline{x2}{x3}\ncline{x3}{x4}\ncline{x4}{x5}\ncline{x5}{x6}
\ncline{x6}{x7}\ncline{x7}{x8}\ncline{x8}{x9}\ncline{x9}{x10}\ncline{x10}{x11}\ncline{x11}{x12}
\ncline{x1}{x12}\ncline{x12}{x2}\ncline{x2}{x4}\ncline{x4}{x6}\ncline{x6}{x8}\ncline{x8}{x10}\ncline{x12}{x10}

\rput(-3.73,4.85){$v_{_3}$}\rput(-2.865,4.35){$u_{_3}$}\rput(-2,3.85){$v_{_4}$}\rput(-1.65,2.5){$u_{_4}$}\rput(-2,1.15){$v_{_5}$}
\rput(-2.865,0.65){$u_{_5}$}\rput(-3.73,0.15){$v_{_0}$}\rput(-4.595,0.65){$u_{_0}$}\rput(-5.46,1.15){$v_{_1}$}
\rput(-5.81,2.5){$u_{_1}$}\rput(-5.46,3.85){$v_{_2}$}\rput(-4.595,4.35){$u_{_2}$}
\rput(-3.73,-0.5){$M(C_6)$}

\cnode(2.73,4.5){3pt}{x1}\cnode(3.595,4){3pt}{x2}\cnode(4.46,3.5){3pt}{x3}\cnode(4.46,2.5){3pt}{x4}\cnode(4.46,1.5){3pt}{x5}
\cnode(3.595,1){3pt}{x6}\cnode(2.73,0.5){3pt}{x7}\cnode(1.865,1){3pt}{x8}\cnode(1,1.5){3pt}{x9}\cnode(1,2.5){3pt}{x10}
\cnode(1,3.5){3pt}{x11}\cnode(1.865,4){3pt}{x12}
\ncline{x1}{x2}\ncline{x3}{x4}\ncline{x5}{x6}
\ncline{x6}{x7}\ncline{x7}{x8}\ncline{x8}{x9}\ncline{x10}{x11}
\ncline{x1}{x12}\ncline{x2}{x4}\ncline{x4}{x6}\ncline{x8}{x10}\ncline{x12}{x10}

\rput(2.73,4.85){$v_{_3}$}\rput(3.595,4.35){$u_{_3}$}\rput(4.46,3.85){$v_{_4}$}\rput(4.81,2.5){$u_{_4}$}\rput(4.46,1.15){$v_{_5}$}
\rput(3.595,0.65){$u_{_5}$}\rput(2.73,0.15){$v_{_0}$}\rput(1.865,0.65){$u_{_0}$}\rput(1,1.15){$v_{_1}$}
\rput(0.65,2.5){$u_{_1}$}\rput(1,3.85){$v_{_2}$}\rput(1.865,4.35){$u_{_2}$}
\rput(2.73,-0.5){$M^{\ast}(C_6)$}
\end{pspicture}
\caption{\small The graphs $M(C_6)$ and $M^{\ast}(C_6)$. \label{fig0}}
\end{center}
\end{figure}

Throughout this section, we will use the following notations (see Fig.\ref{fig0}).

Let $C_{2n}=v_0v_1\cdots v_{2n-1}v_0$. The middle graph of $C_{2n}$, denoted by $M(C_{2n})$,
is obtained from $C_{2n}$ by inserting $u_i$ into $v_iv_{(i+1)\mod(2n)}$,
and connecting $u_iu_{(i+1)\mod(2n)}$ $(0\leq i \leq2n-1)$.
The graph $M^{\ast}(C_{2n})$ is obtained from $M(C_{2n})$ by removing the edges $v_iu_i$ for $1\leq i\leq n-1$,
$u_{n-1}u_n$ and $u_jv_{j+1}$ for $n\leq j\leq 2n-2$.

We define the vertex subsets $A$ and $B$ of $V(M^{\ast}(C_{2n}))$ by
$$A=\{v_1,v_2,\ldots,v_{n-1},u_0,u_1,\ldots,u_{n-1}\},$$ $$B=\{v_{n+1},v_{n+2},\ldots,v_{2n-1},u_n,u_{n+1},\ldots,u_{2n-1}\}.$$

For simplicity, among $M(C_{2n})\times G$ (or $M^*(C_{2n})\times G$), let $p_i=p(v_i\times G)$, $r_i=r(v_i\times G)$, $q_i=p(u_i\times G)$, $s_i=r(u_i\times G)$, $p(A)=p(A\times G)$, $p(B)=p(B\times G)$.

\begin{lem}{\rm{(\cite{m92})}}\label{nl00}
Trees satisfy the two-pebbling property.
\end{lem}
\begin{lem}{\rm{(\cite{yzz})}}\label{lem2.3}
If $n\geq2$, then $f(M(C_{2n}))=2^{n+1}+2n-2$.
\end{lem}

From Lemma 1.2 and the  proof of Lemma 3.2, it is not hard to obtain the following
\begin{lem}\label{lem2.3(1)}
If $n\geq2$, then $f(M^*(C_{2n}),v_0)=2^{n+1}+2n-2$.
\end{lem}
\begin{pps}\label{p1}
$M(C_{2n})$ satisfies the two-pebbling property.
\end{pps}

\begin{pf} Since symmetry, it is clear that
$$f(M(C_{2n}))=\max\{f(M(C_{2n}),v_0), f(M(C_{2n}),u_0)\}.$$
Assume that the target vertex is $v_0$, and $p+k\geq2f(M(C_{2n}))+1$. Since $k\leq4n\leq f(M(C_{2n}))$,
we have $p\geq f(M(C_{2n}))+1$.
Thus if there is one pebble  located on $v_0$, then with the remaining $f(M(C_{2n}))$ pebbles, a
second pebble can be moved to $v_0$.\\

\indent Now,  suppose that $p(v_0)=0$.
We will prove that with the same arrangement of pebbles on $M^*(C_{2n})$, two pebbles can be moved to $v_0$.

Let $H=M^*(C_{2n})$, $C=H[A\setminus v_1]$, and $D=H[B\setminus v_{2n-1}]$.
Then by Lemma~\ref{lem1}, $$f(C)=f(D)=2^{n-1}+n-2,$$ $$f(C\cup \{v_0\})=f(D\cup \{v_0\})=2^n+n-2,$$
$$f(C\cup \{v_n\})=f(D\cup \{v_n\})=2^n+n-2.$$

We consider the worst case, that is $p(v_1)=k(v_1)=p(v_{2n-1})=k(v_{2n-1})=1$ (where $k(v_i)=1$ if there is at least one pebble located on $v_i$ and $0$ otherwise), then
$$p(C)+k(C)+p(D)+k(D)+p_n+q_n+4\geq2^{n+2}+4n-3.$$

If $p(C)+k(C)> 2^{n+1}+2n-4$, then by Lemma~\ref{nl00}, two pebbles can be moved to $v_0$.
Thus we may assume that $p(C)+k(C)\leq 2^{n+1}+2n-4$ and $p(D)+k(D)\leq 2^{n+1}+2n-4$.  We will show that
every one of the vertices $u_0$ and  $u_{2n-1}$ will get at least two pebbles by a sequece of pebbling moves.

Let $p'_n=2^{n+1}+2n-4-p(C)-k(C)\geq 0$, and paint all the pebbles on $C$ red along
with  the $p'_n$ pebbles on $v_n$. Similarly, paint the pebbles on $D$ black,
along with $p''_n=2^{n+1}+2n-4-p(D)-k(D)$ pebbles on $v_n$.
It is easy to see there are enough pebbles on $v_n$ to do this.

Now either $p(C)+k(C)=2^{n+1}+2n-4$ or there are red pebbles on $v_n$.
If equality holds, then $p(C)\geq2^n+n-2$, then two red pebbles can be moved to $u_0$.
If there are red pebbles on $v_n$, then $k'_n=1$, and the red pebbles satisfy
$$p(C)+k(C)+p'_n+k'_n=2^{n+1}+2n-3,$$
and again two red pebbles can be moved to $u_0$.
Similarly, two black pebbles can be moved to $u_{2n-1}$,
so we can move one red pebble and one black pebble to $v_0$.

If the target vertex is $u_0$, then a similar argument can show that there are at least two pebbles can be moved to $u_0$.
\end{pf}


\begin{lem}{\rm{(\cite{h03})}}\label{lem2.01}
Let $P_k=x_1x_2\ldots x_k$ be a path, and let $g$ be some vertex in a graph $G$ which satisfies the odd two-pebbling property. Then, from any arrangement of $(2^k-1)f(G)$ pebbles on $P_k\times G$, it is possible to put a pebble on every $(x_i,g)$ at once $(1\leq i \leq k)$.
\end{lem}

\begin{lem}\label{lem2.1}
Let $P_k=x_1x_2\ldots x_k$ be a path $(k\geq 2)$, and $g$ be some vertex in a graph $G$ which satisfies the odd two-pebbling property. Then from any arrangement of $(2^k-2)f(G)$ pebbles on $x_k\times G$, it is possible to put a pebble on every $(x_i,g)$ at once $(1\leq i \leq k-1)$.\end{lem}

\begin{pf}
We use induction on $k$, where the case $k=2$ is trivial. If it is true for $k-1$, suppose there are $(2^k-2)f(G)$ pebbles on $x_k\times G $,
we use $(2^{k-1}-2)f(G)$ pebbles to put a pebble on every $(x_i,g)$ at once $(2\leq i \leq k-1)$, and with the remaining $2^{k-1}f(G)$
pebbles we can put one pebble on $(x_1,g)$.
\end{pf}

\begin{lem}\label{lem2.2}
Let $T_k$ be the graph obtained from $P_k$ by joining $x_i$ to a new vertex $y_i$ $(1\leq i\leq k-1)$, where $P_k=x_1x_2\ldots x_k$ is a path $(k\geq 2)$.
Let $g$ be some vertex in a graph $G$ which satisfies the odd two-pebbling property.
Then for any arrangement of $(2^k+k-3)f(G)$ pebbles on $T_k\times G$, one of the following will occur\\
\indent  (1)\; we can put a pebble on every $(x_i,g)$ at once $(1\leq i\leq k-1)$;\\
\indent  (2)\; we can put two pebbles on $(x_1,g)$.
\end{lem}

\begin{figure}[ht]
\begin{center}
\hspace*{30pt}
\begin{pspicture}(-3,0)(3,2.5)

\cnode(-2.5,2){3pt}{x1}\cnode(-1.5,2){3pt}{x2}\cnode(0,2){3pt}{x3}\cnode(1,2){3pt}{x4}\cnode(2,2){3pt}{x5}
\cnode(-2.5,1){3pt}{y1}\cnode(-1.5,1){3pt}{y2}\cnode(0,1){3pt}{y3}\cnode(1,1){3pt}{y4}
\ncline{x1}{x2}\ncline{x3}{x4}\ncline{x4}{x5}
\ncline{x1}{y1}\ncline{x2}{y2}\ncline{x3}{y3}\ncline{x4}{y4}
\rput(-.75,2){$\bm{\cdots}$}

\rput(-2.5,2.35){$x_{_1}$}\rput(-1.5,2.35){$x_{_2}$}\rput(0,2.35){$x_{_{k-2}}$}\rput(1,2.35){$x_{_{k-1}}$}\rput(2,2.35){$x_{_{k}}$}
\rput(-2.5,0.65){$y_{_1}$}\rput(-1.5,0.65){$y_{_2}$}\rput(0,0.65){$y_{_{k-2}}$}\rput(1,0.65){$y_{_{k-1}}$}

\end{pspicture}
\caption{\small The graph $T_k$ in Lemma~\ref{lem2.2}. \label{fig1}}
\end{center}
\end{figure}
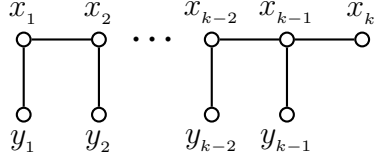
\begin{pf}
While $k=2$, by Lemma~\ref{lem1} and Lemma~\ref{lem2.0}, with $3f(G)$ pebbles on $T_2\times G$,
one pebble can be moved to the vertex $(x_1,g)$.\\
\indent Suppose that there are $(2^k+k-3)f(G)$ pebbles on $T_k\times G$ for $k \geq 3$.
Let $T_{k}'=T_k\setminus \{x_1,y_1\}$. Clearly, $T_k'\cong T_{k-1}$.\\
\indent If $p(T_k'\times G)< (2^{k-1}+k-4)f(G)$,
then $p(P_{x_1y_1}\times G)\geq (2^{k-1}+1)f(G)\geq 5f(G)$. Then clearly, we can move two pebbles to $(x_1,g)$.\\
\indent If $p(T_k'\times G)\geq (2^{k-1}+k-4)f(G)$, and $p_k\geq (2^{k-1}-2)f(G)$. By Lemma~\ref{lem2.1}, using $(2^{k-1}-2)f(G)$ pebbles on $x_k\times G$, we can put a pebble on every $(x_i,g)$ for $2\leq i\leq k-1$.
With the remaining $(2^{k-1}+k-1)f(G)$ pebbles, by Lemma~\ref{lem1} and Lemma~\ref{lem2.0}, we can put one pebble on $(x_1,g)$ for $f(T_k\times G,(x_1,g))\leq f(T_k,x_1)f(G)=(2^{k-1}+k-1)f(G)$.

If $p(T_k'\times G)\geq (2^{k-1}+k-4)f(G)$, and $p_k< (2^{k-1}-2)f(G)$. We use induction in this case,
while the case $k=2$ holds.

Let $r_y$ be the number of vertices with an odd number of pebbles in $\{y_2,y_3,\ldots,y_{k-1}\}\times G$.
We only need to take off $r_y$ pebbles from $\{y_2,y_3,\ldots,y_{k-1}\}\times G$ so that each vertex in it has an even number of pebbles.
It is clear that $r_y\leq (k-2)|V(G)|\leq (k-2)f(G)$, so $r_y+p_k<(2^{k-1}+k-4)f(G)$. So we can choose $(2^{k-1}+k-4)f(G)$ pebbles from $T_k'\times G$ which contains all
pebbles on $x_k\times G$, so that the number of the remaining pebbles on each vertex of $\{y_2,y_3,\ldots,y_{k-1}\}\times G$ is even except at most one vertex.
By induction, with these $(2^{k-1}+k-4)f(G)$ pebbles we can put one pebble on every $(x_i,g)$ at once for $2\leq i\leq k-1$ or move two pebbles to $(x_2,g)$ and then at least one pebble can be moved to $(x_1,g)$.

Now we  prove that with the remaining $(2^{k-1}+1)f(G)$ pebbles, one pebble can be moved to $(x_1,g)$.

Let $\tilde{p}_y=\sum\limits_{i=2}^{k-1}\tilde{p}(y_i\times G)$. Let $P_1$ denote the path $y_1x_1x_2\ldots x_{k-1}$, and
$P_2$ denote the path $x_1x_2\ldots x_{k-1}$.

For the number of the remaining pebbles on each vertex of $\{y_2,y_3,\ldots,y_{k-1}\}\times G$ is even except at most one vertex, then we can move $\left\lfloor\frac{1}{2}\tilde{p}_y\right\rfloor$ pebbles from the vertices of $\{y_2,y_3,\ldots,y_{k-1}\}\times G$ to $\{x_2,x_3,\ldots,x_{k-1}\}\times G$.

\textbf{Case $1$.}\, $\tilde{p}_y\leq 2^{k-1}f(G)-1$.
Then $$\tilde{p}(P_1\times G)=(2^{k-1}+1)f(G)-\tilde{p}_y+\left\lfloor\frac{1}{2}\tilde{p}_y\right\rfloor\geq(2^{k-2}+1)f(G).$$
By Lemma~\ref{lem2.0}, $f(P_1\times G,(x_1,g))\leq f(P_1,x_1)f(G)=(2^{k-2}+1)f(G)$, so one pebble can be moved to $(x_1,g)$.

\textbf{Case $2$.}\, $\tilde{p}_y\geq 2^{k-1}f(G)$. Then
$\tilde{p}(P_2\times G)\geq \left\lfloor\frac{1}{2}\tilde{p}_y\right\rfloor\geq 2^{k-2}f(G).$
By Lemma~\ref{lem2.0}, $f(P_2\times G,(x_1,g))\leq f(P_2,x_1)f(G)=2^{k-2}f(G)$, so one pebble can be moved to $(x_1,g)$.
\end{pf}

\begin{thm}\label{thm1}
If $G$ satisfies the odd two-pebbling property, then
$$f(M(C_{2n})\times G)\leq f(M(C_{2n}))f(G)=(2^{n+1}+2n-2)f(G).$$
\end{thm}

\begin{pf}
Suppose that there are  $(2^{n+1}+2n-2)f(G)$ pebbles placed on the vertices of $M(C_{2n})\times G$, we will show that at least one pebble can be moved to the target vertex.\\
\indent Since symmetry, it is clear that
$$f(M(C_{2n})\times G)=\max\{f(M(C_{2n})\times G,(v_0,g)), f(M(C_{2n})\times G,(u_0,g))\}.$$ So we only need to distinguish  two cases.\\
\indent \textbf{Case $1$.} The target vertex is $(v_0,g)$.\\
\indent\textbf{Subase $1.1$.}  $p_n+r_n\leq 2f(G)$.\\
\indent  We remove all the pebbles off $v_n\times G$ such that
\begin{align*}
\tilde{p}((M^\ast(C_{2n})\setminus v_n)\times G)&=\frac{p_n-r_n}{2}+p(A)+p(B)+p_0\\
&=-\frac{1}{2}(p_n+r_n)+p_n+p(A)+p(B)+p_0\\
&\geq  (2^{n+1}+2n-3)f(G).\end{align*}
By Lemma~\ref{lem1}, $f(M^\ast(C_{2n})\setminus v_n,v_0)=2^{n+1}+2n-3$. According to Lemma~\ref{lem2.0}, one pebble can be moved to $(v_0,g)$.\\
\indent\textbf{Subase $1.2$.} $p_n+r_n>2f(G)$.\\
\indent  Then we can put two pebbles to $(v_n,g)$. Note that $p_n$ and $r_n$ are of the same-parity,  we keep $2f(G)-r_n+2$
 pebbles on $v_n\times G$ so that at least two pebbles still can be moved to $(v_n,g)$, and move the rest pebbles to $A\times G$.
 So
$$\tilde{p}(A\times G)=\frac{1}{2}(p_n-(2f(G)-r_n+2))+p(A)=\frac{p_n+r_n}{2}-f(G)-1+p(A). \eqno(3.1)$$

By Lemma~\ref{lem1}, $f(M^*(C_{2n})[B,v_0],v_0)=2^n+n-1$, so if we move as many as possible pebbles from $v_n\times G$ to $B\times G$, then one pebble can be moved to $(v_0,g)$ unless
$$\frac{p_n-r_n}{2}+p(B)+p_0\leq (2^n+n-1)f(G)-1.\eqno(3.2)$$

If (3.2) holds, then
$$\begin{array}{ll}
\tilde{p}(A\times G)&=\frac{1}{2}(p_n+r_n)-f(G)-1+p(A)\\
&\geq  p_n+p(A)+p(B)+p_0-f(G)-(2^n+n-1)f(G)\\
&=(2^{n+1}+2n-2)f(G)-(2^n+n)f(G)\\
&=(2^n+n-2)f(G).
\end{array}\eqno(3.3)$$

\noindent It follows from Lemma~\ref{lem1} that
$$f(M^*(C_{2n})\setminus\{v_n,u_{n-1},v_{n-1}\},v_0)=3\cdot2^{n-1}+2n-4.$$
Thus if we move as many as possible pebbles from
$v_n\times G$ to $u_n\times G$, and from $v_{n-1}\times G$ to $u_{n-2}\times G$, then one pebble can be moved to $(v_0,g)$
unless
$$\begin{array}{ll}
\frac{1}{2}(p_n-r_n)+\frac{1}{2}(p_{n-1}-r_{n-1})+p(B)+p_0+(p(A)-p_{n-1}-q_{n-1})\\
\hspace{3cm}\leq (3\cdot2^{n-1}+2n-4)f(G)-1.
\end{array}\eqno(3.4)$$

If (3.4) holds, then
$\frac{1}{2}(p_n+r_n)+\frac{1}{2}(p_{n-1}+r_{n-1})+q_{n-1}\geq (2^{n-1}+2)f(G)+1.$\\
Thus
$$\left(\frac{1}{2}(p_n+r_n)-f(G)-1+q_{n-1}\right)+\left(\frac{1}{2}(p_{n-1}+r_{n-1})-f(G)-1\right)\geq 2^{n-1}f(G)-1.\eqno(3.5)$$
\indent\textbf{Subase $1.2.1$.}\, $\frac{1}{2}(p_n+r_n)-f(G)-1+q_{n-1}\geq f(G)$. \\
\indent Then from (3.3) it follows that, with $f(G)$ pebbles on $u_{n-1}\times G$, one pebble can be moved to $(u_{n-1},g)$; and from Lemma~\ref{lem2.2} it follows that, with the remaining $(2^n+n-3)f(G)$ pebbles, we can put one pebble to each $(u_i,g)$ for $0\leq i\leq n-2$ or put two pebbles to $(u_0,g)$, we can move one more pebble to
$(u_{n-1},g)$ with $2f(G)-r_n+2$ pebbles on $v_n\times G$, so one pebble can be moved to $(v_0,g)$.

\indent\textbf{Subase $1.2.2$.} $\frac{1}{2}(p_n+r_n)-f(G)-1+q_{n-1}< f(G)$.\\
\indent  Then from (3.5), we have
$$\frac{p_{n-1}+r_{n-1}}{2}-f(G)-1\geq (2^{n-1}-1)f(G).\eqno(3.6)$$
So we can keep $2f(G)-r_{n-1}+2$ pebbles on $v_{n-1}\times G$, so that one pebble can be moved to $(u_{n-2},g)$, and moving no less than $(2^{n-1}-1)f(G)$ pebbles to $u_{n-2}\times G$. With these pebbles, by Lemma~\ref{lem2.01}, we can put one pebble to every
$(u_i,g)$ at once $(0\leq i\leq n-2)$. So one pebble can be moved to $(v_0,g)$.\\

\textbf{Case $2$.} The target vertex is $(u_0,g)$.

Let $M'(C_{2n})$ be the graph obtained from $M(C_{2n})$ by removing the edges $u_iv_{i+1}$ for $0\leq i\leq n-2$ and $u_jv_j$
 for $n+2\leq j\leq 2n-1$ and $u_nv_n$, $u_nv_{n+1}$, $u_0v_0$.

Let $A'=\{u_1,u_2,\ldots,u_{n-1},v_1,v_2,\ldots,v_n\}$ and
$B'=\{u_{n+1},u_{n+2},\ldots,u_{2n-1},v_{n+1},v_{n+2},\\
\ldots,v_{2n-1},v_0\}.$\\

It is clear that $M'(C_{2n})[A']\cong M^{\ast}(C_{2n})[A]$,
and $M'(C_{2n})[B',u_0]\cong M^{\ast}(C_{2n})[B,v_0]$.
We only need to prove that
one pebble can be moved from $M'(C_{2n})\times G$ to $(u_0,g)$.\\
\indent\textbf{Subcase $2.1$.} $q_n+s_n\leq 2f(G)$.\, By a similar process as before, one pebble can be moved to $(u_0,g)$.\\
\indent\textbf{Subcase $2.2$.} $q_n+s_n> 2f(G)$.\\
\indent Then by a similar process as before, we can keep $2f(G)-s_n+2$
 pebbles on $u_n\times G$ so that two pebbles can be moved to $(u_n,g)$, and move the rest pebbles to $A'\times G$.
 So
 \begin{align*}
 \tilde{p}(A'\times G)=&\frac{q_n-(2f(G)-s_n+2)}{2}+p(A')\\
 =&\frac{q_n+s_n}{2}-f(G)-1+p(A').
\end{align*}

Similarly, if we move as many as possible pebbles from $u_n\times G$ to $B'\times G$, then one pebble can be moved from
$B'\times G$ to $(u_0,g)$, unless
\begin{align*}
 \tilde{p}(A'\times G)\geq& (2^n+n-2)f(G).
\end{align*}

According to Lemma~\ref{lem2.2}, we can put one pebble on $(u_i,g)$ at once for
$1\leq i\leq n-1$ or put two pebbles on $(u_1,g)$. With $2f(G)-s_n+2$ pebbles on $u_n\times G$, one more pebble can be moved to $(u_{n-1},g)$. So one pebble can be moved to $(u_0,g)$.
\end{pf}

\section{Remark}
In this paper, we show that $f(C_{4k+3}\times G)\leq f(C_{4k+3})f(G)$, where $G$ satisfies the odd two-pebbling property.
But we cannot show $f(C_{4k+1}\times G)\leq f(C_{4k+1})f(G)$. Maybe since the odd cycles have not been solved completely, the middle graph of the odd cycle is an open problem, too.

\end{document}